\documentclass[12pt,epsfig,epsf]{article}
\usepackage{url}
\usepackage{epic, eepic, subfigure, floatflt}
\usepackage{amsfonts, epsfig, latexsym, amsmath}
\begin{document}
\newcommand{\R}{\ensuremath{\mathbb{R}}}
\newcommand{\N}{\ensuremath{\mathbb{N}}}
\newcommand{\Q}{\ensuremath{\mathbb{Q}}}
\newcommand{\C}{\ensuremath{\mathbb{C}}}
\newcommand{\Z}{\ensuremath{\mathbb{Z}}}
\newcommand{\T}{\ensuremath{\mathbb{T}}}
\newtheorem{theorem}{Theorem}[section]
\newtheorem{definition}[theorem]{Definition}
\newtheorem{conjecture}[theorem]{Conjecture}
\newtheorem{corollary}[theorem]{Corollary}
\newtheorem{lemma}[theorem]{Lemma}
\newtheorem{claim}[theorem]{Claim}
\newtheorem{remark}[theorem]{Remark}
\newtheorem{proposition}[theorem]{Proposition}
\newcommand{\qed}{\hfill $\Box$ }
\newcommand{\proof}{\noindent{\bf Proof.}\ \ }
\newcommand{\sketchproof}{\noindent{\bf Sketch of proof.}\ \ }
\def\QuotS#1#2{\leavevmode\kern-.0em\raise.2ex\hbox{$#1$}\kern-.1em/\kern-.1em\lower.25ex\hbox{$#2$}}

\title{\Large{\bf Zigzag structure of complexes}}

\author{Michel DEZA\\
        \normalsize   LIGA, ENS, Paris and Institute of Statistical Mathematics, Tokyo\\
\and
        Mathieu DUTOUR\\
        \normalsize  LIGA, ENS, Paris and Hebrew University, Jerusalem
\footnote{Research financed by EC's IHRP Programme, within the Research Training Network ``Algebraic Combinatorics in Europe,'' grant HPRN-CT-2001-00272.}\\
}


\date{\small \today}
 
\maketitle

\parindent=1cm

\begin{abstract}
Inspired by Coxeter's notion of Petrie polygon for $d$-polytopes (see \cite{Cox73}),
we consider a generalization of the notion of zigzag circuits on complexes
and compute the zigzag structure for several interesting families of
$d$-polytopes, including semiregular, regular-faced, Wythoff Archimedean
ones, Conway's $4$-polytopes, half-cubes, folded cubes. 

Also considered are regular maps and Lins triality relations on maps.
\end{abstract}

{\small {\em Mathematics Subject Classification}: Primary 52B05, 52B10;
Secondary 05C30,\\
\indent 05C10.

}

\bigskip
\bigskip

\section{Introduction}
The notion of zigzag was introduced for plane graphs in \cite{Sh} (as a {\em left-right} path) and for regular polytopes in \cite{Cox73} (as a {\em Petrie polygon}). 
We focus here on generalization of zigzags for higher dimension.

Zigzags can be also defined for maps on orientable surface;
see, for example, on Figure \ref{Dyck-Klein-Map} typical zigzags for
dual Klein map $\{7,3\}$ and dual Dyck map
$\{8,3\}$.
Moreover, this notion, being local, is defined even for non-oriented maps.
See Section \ref{LinsSection} on maps.
Also, the notion of zigzag extends naturally on infinite plane graphs.

We use for polytopes notations and terminology of \cite{Cox73}; for example, $\alpha_d$, $\beta_d$, $\gamma_d$ and $\frac{1}{2}\gamma_d$ denote $d$-dimensional simplex, cross-polytope, cube and half-cube, respectively.
Their $1$-skeleton graphs are denoted by $K_{d+1}$, $K_{d\times 2}$, $H_d$ and $\frac{1}{2}H_d$, respectively.
We use also Schl\"afli notation from \cite{Cox73} in Tables \ref{RegularResult} and \ref{CoxeterTable}.
By $Prism_m$ and $APrism_m$ are denoted semiregular $m$-gonal prism and $m$-gonal antiprism, respectively.

The {\em medial} of a polytope $P$, denoted by $Med(P)$, is the polytope formed by the convex hull of the midpoints of all edges of $P$.
It can also be defined combinatorially on maps on surfaces by taking as vertices the edge of the original map, by taking as edges the pair of edges sharing an incident vertex and an incident face and by taking as faces the vertices and faces of the original map.
This notion of medial can also be defined combinatorially on $d$-dimensional complexes, including {\em maps}, i.e. the case $d=2$.





\section{Zigzags for $d$-dimensional complexes}\label{GeneralCase}
We extend here the definitions of zigzags to any complex.
A {\em chain} of length $k$ in a partially ordered set is a sequence $(x_0,\dots x_k)$, such that $x_{i} < x_{i+1}$.
A chain $C$ is a {\em subchain} of another chain $C'$ if it is obtained by removing 
some elements in $C'$.

A chain is {\em maximal} if it is not a subchain of another chain.
The {\em rank} $rank(x)$ of an element $x$ is the maximal length of chains, beginning at the lowest elements $0$ and terminating at $x$.
A partially ordered set is called {\em ranked} if there is a lowest element $0$ and a greatest element $1$ and if,
given two elements $x<y$ with no elements $z$ satisfying to $x<z<y$, one has $rank(y)=1+rank(x)$.

A partially ordered set is called a {\em lattice} if for any
two elements $x$ and $y$, there are an unique smallest element $s$
and an unique greatest element $t$, such that $x\leq s$, $y\leq s$ 
and $x\geq t$, $y\geq t$.

The {\em dimension} of an element is defined as $rank(x)-1$.

A {\em $d$-dimensional complex} ${\cal K}$ is a finite partially ordered set, such that it holds:

(i) ${\cal K}$ has a smallest $0$ and highest element $1$,

(ii) ${\cal K}$ is ranked and all maximal chains have length $d+2$,

(iii) given two elements $x$ and $y$ with $x\leq y$ and $dim(y)=2+dim(x)$,
there are exactly two elements $u, u'$, such that $x\leq u\leq y$ and
$x\leq u'\leq y$.

A $d$-dimensional complex is called {\em simplicial} if for every element
$x$ of dimension $d$, there is exactly $d+1$ elements of dimension $0$
contained in it.

In a $d$-dimensional complex, a maximal chain is called a {\em flag}; it
necessarily begins at $0$ and terminates at $1$.

Using (iii), one can define the following permutation operator on flags.
For $1\leq i\leq d+1$, denote by $\sigma_i$ the operator transforming
$(0,x_1,\dots,x_i,\dots, x_{d+1},1)$ into the flag
$(0,x_1,\dots,x'_i,\dots, x_{d+1},1)$ with $x'_i$
being the unique element satisfying to $x'_i\not= x_i$
and $x_{i-1}\leq x'_i\leq x_{i+1}$.
One has $\sigma_i^2=1$ and $\sigma_i\sigma_j=\sigma_j\sigma_i$ if $i<j-1$.

\begin{definition}
Let ${\cal K}$ be a $d$-dimensional complex, then:

(i) denote by ${\cal F}({\cal K})$ the set of flags of ${\cal K}$,

(ii) denote by ${\cal G}({\cal K})$ the graph having, as vertex-set,
${\cal F}({\cal K})$, with two flags being adjacent if they are obtained one 
from the other by a permutation $\sigma_i$,

(iii) the complex ${\cal K}$ is said to be {\em orientable} if
${\cal G}({\cal K})$ is bipartite; an {\em orientation} of ${\cal K}$
consists in selecting one of the two connected components.
\end{definition}
In the case $d=2$, the elements of dimension $0$, $1$ and $2$
are called {\em vertices}, {\em edges} and {\em faces}, respectively.

The definition of orientability, given above, corresponds to the fact
that, given an orientation on a cell complex and a maximal chain
$(f_1,\dots, f_{d-1})$ of faces, one can find the last face $f_d$
that makes it a flag.

A {\em $(d+1)$-polytope} $P$ is defined as the convex hull of a set of points in $\R^{d+1}$. The set of faces of $P$ defines a lattice and so, a $d$-dimensional complex, which is a lattice, since the boundary of a $d$-polytope is homeomorphic to $S^{d-1}$.



Call a $d$-dimensional complex {\em regular} if its symmetry group is
transitive on the set of flags.

\begin{theorem}\label{FundamentalTheorem}
Let ${\cal K}$ be a $d$-dimensional complex and $x=(0,x_1,\dots,x_{d+1},1)$ be a flag in ${\cal K}$.

Then there exists an unique sequence of faces $(x_{i,j})_{1\leq i\leq d+1, 1\leq j\leq d+2-i}$, namely:
\begin{equation*}
\begin{array}{c}
x_{1,1},\hspace{1mm}\dots\dots, \hspace{1mm}x_{1, j}, \hspace{1mm}\dots\dots, \hspace{1mm}x_{1,d+1}\\
x_{2,1}, \hspace{1mm}\dots\dots, \hspace{1mm}x_{2, j}, \hspace{1mm}\dots\dots, \hspace{1mm}x_{2,d}\\
\vdots\\
x_{d-1,1},x_{d-1,2},x_{d-1,3}\\
x_{d,1},x_{d,2}\\
x_{d+1,1}
\end{array}
\end{equation*}
such that it holds:

(i) $x_{i,1}=x_i$,

(ii) $dim(x_{i,j})=i-1$,

(iii) $x_{i,j}\leq x_{i+1,j}$ for $1\leq i\leq d$ and $1\leq j\leq d+1-i$,

(iv) $x_{i,j}\leq x_{i+1,j-1}$ for $1\leq i\leq d$ and $2\leq j\leq d+2-i$.

Moreover, if ${\cal K}$ is a lattice, then the elements $(x_{i,j})$ are uniquely defined by the vertex sequence $(x_{1,j})_{1\leq j\leq d+1}$.

\end{theorem}
\proof Using property (iii), one can find successively, $x_{1,2}$, \dots, $x_{d,2}$, then $x_{1,3}$ and so on.

If ${\cal K}$ is a lattice, then $x_{i,j}$ can be characterized as 
the smallest element greater than $x_{i-1,j}$ and $x_{i-1,j+1}$. \qed

\begin{definition}
Let ${\cal K}$ be a $d$-dimensional complex.

(i) Denote by $T=\sigma_{d+1}\sigma_{d}\dots \sigma_1$ the {\em translation operator} of ${\cal K}$.

(ii) A {\em zigzag} in ${\cal K}$ is a circuit $(f_1, \dots, f_l)$ of flags, such that
$f_{j+1}=T(f_j)$; $l$ denotes the length of the zigzag.

(iii) Given a flag $f$, the {\em reverse} $f^{t}$ of $f$ is defined as $(0$, $x_{1,d+1}$, $x_{2,d}$,\dots, $x_{d,2}$, $x_{d+1, 1}$, $1)$ with $(x_{i,j})$ as in Theorem \ref{FundamentalTheorem}.

(iv) The {\em reverse} of a zigzag $(f_1,\dots, f_l)$ is the zigzag $(f_{l}^{t}, f_{l-1}^t,\dots, f_1^{t})$.
\end{definition}
The above notion (central in this paper), for the special case of
an $d$-polytope,
essentially coincides with the following notion on page 223 of \cite{Cox73}: ``A Petrie polygon of an $d$-dimensional polytope or of an $(d-1)$-dimensional honeycomb, is a skew polygon, such that any $(d-1)$ consecutive sides but no $d$, belong to a Petrie polygon of a cell.''

The choice of a zigzag $(f_1, \dots, f_l)$ over its reverse $(f_l^{t}, \dots, f_1^{t})$ amounts to choosing an orientation on the zigzag.
In the sequel a zigzag is identified with its reverse.


Note that if ${\cal K}$ is a $d$-dimensional simplicial complex with $f$ facets, then one has $|{\cal F}({\cal K})|=(d+1)!f$.
Note also that the stabilizer of a flag is trivial and so, if ${\cal K}$ has $p$ orbits of flags, then $|{\cal F}({\cal K})|=p|Sym({\cal K})|$.


\begin{proposition}
If the complex ${\cal K}$ is oriented and of even dimension, then the length of any zigzag is even.
\end{proposition}
\proof Since ${\cal K}$ is oriented, the set ${\cal F}({\cal K})$ is splitted in two parts, ${\cal F}_1$ and
${\cal F}_2$.
Since $d$ is even, the translation $T=\sigma_{d+1}\sigma_{d}\dots \sigma_1$
of all its flags has an odd number of components; so, it interchanges ${\cal F}_1$ and ${\cal F}_2$. \qed

On Section 5.91 of \cite{Cox73} the evenness of the length of zigzags was obtained for complexes arising from Coxeter groups of dimension $3$;
there was given the formula $g=h(h+2)$ with $g$ being the size of the group and $h$ the length of the zigzag.

\begin{definition}
Take a zigzag $Z=(f_1, \dots, f_l)$ and its reverse $Z^t=(f^{t}_l, \dots, f_1^t)$.

(i) Given a flag $f_j$, 
if $\sigma_1(f_j)$ belongs to $Z$, then self-intersection is called of {\em type I}, while if $\sigma_1(f_j)$ belongs to $Z^t$, then it is called of {\em type II}.

(ii) The {\em signature} of the zigzag $Z$ is the pair $(n_I, n_{II})$
with $n_I$ being the number of self-intersections of type I and $n_{II}$
the number of self-intersections of type II. The signature does not change
if one interchanges $Z$ and $Z^t$.

(iii) Take two zigzags $Z_1$ and $Z_2$ with associated circuits
$(f_{1,1}$,\dots, $f_{1,l})$, $(f^{t}_{1,l}$,\dots, $f^t_{1,1})$ and
$(f_{2,1}$,\dots, $f_{2,l})$, $(f^{t}_{2,l}$,\dots, $f^t_{2,1})$. 
If $f'_{1,j}$ belongs to $Z_2$, then it is called an intersection of {\em type I}, while if it belongs to $Z^t_{2}$, it is called an intersection of {\em type II}.

(iv) The {\em signature} $(n_I, n_{II})$ is the pair enumerating such intersections. If $Z_2$ and $Z_{2}^t$ are interchanged, then the types of intersections are interchanged also.

\end{definition}

The {\em $z$-vector} of a complex ${\cal K}$ is the vector enumerating the 
lengths of all its zigzags with their signature as subscript.
The simple zigzags are put in the beginning, in increasing order 
of length, without their signature $(0,0)$, and separated by 
a semicolon from others. Self-intersecting zigzags are also ordered by
increasing lengths. If there are $m>1$ zigzags of the same
length $l$ and the same signature $(\alpha_1, \alpha_2)\not= (0,0)$, 
then we write $l_{\alpha_1,\alpha_2}^m$.
It turns out, that Snub Cube, Snub Dodecahedron, $Pyr(\beta_{d-1})$
and $BPyr(\alpha_{d-1})$ are the only polytopes in Tables of this
paper, having self-intersecting zigzags.

Given two zigzags $Z$ and $Z'$, their {\em normalized signature} is the pair $(n_I, n_{II})$ enumerating intersection of type I and II with orientation chosen so that $n_I\leq n_{II}$.
For a zigzag $Z$, its {\em intersection vector} $Int(Z)=\dots,(c_{k,I}, c_{k,II})^{m_k},\dots$ is such that $(\ldots,(c_{k,I},c_{k,II}),\ldots)$ is a sequence $(c_{k,I},c_{k,II})$ of its non-zero normalized signature with all others zigzags, and $m_k$ denote respective multiplicities. If the zigzag has signature $(n_1, n_{II})$, then its length $l$ satisfies to
\begin{equation*}
l=2(n_I+n_{II})+\sum_{k} m_k(c_{k,I}+c_{k,II})\;.
\end{equation*}

The {\em dual} ${\cal K}^*$ of a $d$-dimensional complex ${\cal K}$ is the complex with the same elements as ${\cal K}$, but with $x\leq y$ in ${\cal K}^*$ being equivalent to $y\leq x$ in ${\cal K}$.

\begin{theorem}
Every zigzag $Z$ in ${\cal K}$ corresponds to an unique zigzag $Z^*$ in ${\cal K}^*$ with the same length.

\end{theorem}
\proof Given a flag $f=(0,x_1,\dots, x_{d+1}, 1)$ of ${\cal K}$,
one can associate to it a flag $f'=(1, x_{d+1},\dots, x_1, 0)$ of ${\cal K}^*$.
Denote by $\sigma'_i$ the operator on ${\cal K}^*$, which acts by changing
the $i$-th element. It is easy to see that its action on $f'$ corresponds
to the action of $\sigma_{d+2-i}$ on $f$. So, one has $T'(f')=(T^{-1}f)'$
and every zigzag $(f_1, \dots, f_l)$ of ${\cal K}$ corresponds
to a zigzag $(f'_l,\dots, f'_1)$ of ${\cal K}^*$. \qed

In the case of maps (i.e. for $d=2$), every intersection in ${\cal K}$ corresponds to an intersection in ${\cal K}^*$ with type I or II interchanged.
This is not, a priori, the case of complexes of dimension $d>2$.

A $d$-dimensional complex ${\cal K}$ is said to be {\em $z$-transitive} if its symmetry group $Sym({\cal K})$ is transitive on zigzags.
It is said to be {\em $z$-knotted} if it has only one zigzag, Note that the stabilizer of a flag is necessarily the trivial group, i.e., every orbit of flags has the size $|Sym({\cal K})|$.

Denote by $Z({\cal K})$ the graph formed by the set of zigzags of 
a complex ${\cal K}$ with two zigzags being adjacent if the signature
of their intersection is different from $(0,0)$.
In the case of a $2$-dimensional complexes, we prove (see Section
\ref{LinsSection}) that $Z({\cal K})$ is connected. In the case of 
complexes of dimension $d>2$, there is no reason to think that 
connectivity will still hold.

\begin{proposition}
If a $d$-dimensional complex ${\cal K}$ is regular, then:

(i) ${\cal K}$ is $z$-transitive,

(ii) if $Z({\cal K})$ is connected, then either zigzags have no self-intersections, or ${\cal K}$ is $z$-knotted.
\end{proposition}
\proof The transitivity on zigzags is obvious. If a zigzag has a self-intersection, then, by transitivity, all flags correspond to a self-intersection of zigzags. Since $Z({\cal K})$ is connected, it means that there is only one zigzag. \qed


\begin{conjecture}
The signature of any zigzag in any odd-dimensional complex is $(0,0)$.
\end{conjecture}
The above conjecture is strange and we do not see why it would be true.
Nevertheless, we did not find a single example violating it.

\section{Some generalizations of regular $d$-polytopes}

Remind, that a {\em regular $d$-polytope} is one whose symmetry group
is transitive on flags.

A {\em regular-faced $d$-polytope} is one having only regular facets.
A {\em semiregular $d$-polytope} is a regular-faced
$d$-polytope whose symmetry group is transitive on vertices.
All semiregular, but not Platonic, $3$-polytopes (i.e. $13$ Archimedean
$3$-polytopes and $Prism_m$, $APrism_m$ for any $m\geq 3$) were discovered
by Kepler (\cite{kepler}).
The list of all $7$ semiregular, but not regular, $d$-polytopes with
$d\geq 4$ was given by Gosset in 1897 (\cite{gosset}), but proofs
were never published; see also \cite{BlBl}. This list consists of $5$
polytopes, denoted by $n_{21}$ (where $n\in \{0,1,2,3,4\}$) of dimension
$n+4$, and two exceptional ones (both $4$-dimensional):
{\em snub $24$-cell} $s(3,4,3)$ and {\em octicosahedric polytope}.
$0_{21}$, $24$-cell, $s(3,4,3)$ and the octicosahedric polytope are
the medials of $\alpha_4$, $\beta_4$, $24$-cell and $600$-cell, respectively
(see also Section \ref{WythoffSection} and Table \ref{WythoffResult} for
the notion of Wythoff Archimedean).
$s(3,4,3)$ is obtained also by eliminating some $24$ vertices of 
$600$-cell (see \cite{Cox73}). $1_{21}$ is $\frac{1}{2}\gamma_5$;
$2_{21}$ and $3_{21}$ are Delaunay polytopes of the root lattices $E_6$
and $E_7$.
The skeleton of $4_{21}$ is the root graph of all $240$ roots of the 
root system $E_{8}$.

\begin{table}
\begin{center}
\begin{tabular}{||c|c|c|c||}
\hline\hline
dimension&complex                      & $z$-vector          &int. vectors\\
\hline
$d-1$&$d$-simplex $\alpha_d$=$\{3^{d-1}\}$         &$(d+1)^{d!/2}$  &$(0,1)^{d+1}$ if $d\geq 4$\\
     &                                             &                &$(1,1)^2$ if $d=3$\\
$d-1$&cross-$d$-polytope=$\beta_d$=$\{3^{d-2},4\}$          &$(2d)^{2^{d-2}(d-1)!}$&$(0,2)^d$\\
$2$&Dodecahedron=$\{5,3\}$       &$10^6$                   &$(0,2)^5$\\
$2$&Great Dodecahedron=$\{5,\frac{5}{2}\}$   & $6^{10}$  &$(0,2)^3$\\
$2$&Petersen graph on $P^2$      &$5^6$                &$(0,1)^5$\\
\hline
$3$&$600$-cell=$\{3,3,5\}$         &$30^{240}$   &$(0,2)^{15}$\\
$3$&$24$-cell=$\{3,4,3\}$          &$12^{48}$&$(0,2)^{6}$\\
\hline
$3$&snub $24$-cell=$s(3,4,3)$      &$20^{144}$   &$(1,1)^{4}, (0,2)^{4}, (0,4)$\\
$3$&octicosahedric polytope        &$45^{480}$  &$(0,1)^{15}, (0,2)^{15}$\\
$3$&$0_{21}$=Med($\alpha_4$)       &$15^{12}$    &$(1,2)^5$\\
$4$&$1_{21}=\frac{1}{2}\gamma_5$=$Med(\beta_5)$  &$12^{240}$ &$(0,1)^{8}, (0,2)^{2}$\\
$5$&$2_{21}$=Schl\"afli polytope (in $E_6$) &$18^{4320}$  &$(0,1)^{6}, (0,2)^{6}$\\
$6$&$3_{21}$=Gosset polytope (in $E_7$)    &$90^{48384}$ &$(0,2)^{15}, (0,4)^{15}$\\
$7$&$4_{21}$ ($240$ roots of $E_8$)                    &$36^{29030400}$  &$(0,1)^{24}, (0,4)^{3}$\\
\hline
$2$&$92$ Johnson solids &\multicolumn{2}{c||}{See Remark \ref{WebPageInformation}}\\
$3$&$Pyr(Icosahedron)$           &$25^{12}$  &$(0,10), (0,3)^{5}$\\
$3$&$BPyr(Icosahedron)$          &$40^{12}$  &$(0,20), (0,4)^{5}$\\
$3$&$0_{21}+Pyr(\beta_3)$        &$42^6$  &$(1,1), (8,8), (12, 12)$\\
$3$&special cuts of $600$-cell&\multicolumn{2}{c||}{See Remark \ref{WebSpecialCut}}\\
$d-1$&$Pyr(\beta_{d-1})$&\multicolumn{2}{c||}{See Conjecture \ref{ConjectureBetaAlpha} }\\
$d-1$&$BPyr(\alpha_{d-1})$&\multicolumn{2}{c||}{See Conjecture \ref{ConjectureBetaAlpha} }\\
\hline
$3$&$45$ Wythoff Archimedean&\multicolumn{2}{c||}{See Table \ref{WythoffResult}}\\
   &$4$-polytopes&\multicolumn{2}{c||}{}\\
$3$&$17$ prisms on Platonic&\multicolumn{2}{c||}{See Table \ref{PrismPlatonicArchimedean}}\\
   &and Archimedean solids&\multicolumn{2}{c||}{}\\
$3$&Grand Antiprism         &\multicolumn{2}{c||}{See Remark \ref{GrandAntiprism}}\\
$3$&$C_p\times C_q$&\multicolumn{2}{c||}{See Conjecture \ref{ProductOfPolygons} }\\
$3$&prisms on $APrism_m$&\multicolumn{2}{c||}{See Conjecture \ref{ConjecturePrismOfAntiprism} }\\
\hline
\hline
\end{tabular}
\end{center}
\caption{$z$-structure of regular, semiregular, regular-faced $d$-polytopes and Conway's $4$-polytopes}
\label{RegularResult}
\end{table}

\begin{table}
{\small
\begin{center}
\begin{tabular}{||c||c|c||c|c||}
\hline\hline
polyhedron $P$                &\multicolumn{2}{|c||}{$P$}      &\multicolumn{2}{|c||}{$Prism(P)$}\\
\hline
                              &$z$         &int. vectors        &$z$        &int. vectors\\
\hline
Tetrahedron                   &$4^3$       &$(1,1)^{2}$         &$16^{6}$   &$(3,3)^{2}, (0,4)$\\
Octahedron                    &$6^4$       &$(0,2)^3$           &$8^{24}$   &$(0,2]^{4}$\\
Dodecahedron                  &$10^6$      &$(0,2)^{5}$         &$40^{12}$  &$(0,6)^{5}, (0,10)$\\
Icosahedron                   &$10^6$      &$(0,2)^{5}$         &$40^{12}$  &$(0,6)^{5}, (0,10)$\\
\hline
Cuboctahedron                 &$8^6$       &$(0,2)^4$           &$32^{12}$  &$(0,6)^{4}, (0,8)$\\
Icosidodecahedron             &$10^{12}$   &$(0,2)^5$           &$40^{24}$  &$(0,6)^{5}, (0,10)$\\
Truncated Tetrahedron         &$12^3$      &$(3,3)^{2}$         &$16^{18}$  &$(0,3)^{4}, (0,4)$\\
                              &            &                    &           &or $(3,3)^{2}, (0,4)$\\
Truncated Octahedron          &$12^6$      &$(0,4), (0,2)^{3}$  &$16^{36}$  &$(0,2)^{4}, (0,4)^{2}$\\
Truncated Cube                &$18^4$      &$(2,4)^{3}$         &$24^{24}$  &$(0,2)^{3}, (0,4)^{3}, (0,6)$\\
                              &            &                    &           &or $(2,4)^3, (0,6)$\\
Truncated Icosahedron         &$18^{10}$   &$(0,2)^{9}$         &$24^{60}$  &$(0,2)^{9}, (0,6)$\\
Truncated Dodecahedron        &$30^6$      &$(2,4)^{5}$         &$40^{36}$  &$(0,2)^{5}, (0,4)^{5}, (0,10)$\\
                              &            &                    &           &or $(2,4)^5, (0,10)$\\
Rhombicuboctahedron           &$12^8$      &$(0,2)^{6}$         &$16^{48}$  &$(0,2)^{6}, (0,4)$\\
Rhombicosidodecahedron        &$20^{12}$   &$(0,2)^{10}$        &$80^{24}$  &$(0,6)^{10}, (0,20)$\\
Truncated Cuboctahedron       &$18^8$      &$(0,2)^{6}, (0,6)$  &$24^{48}$  &$(0,2)^{6}, (0,6)^{2}$\\
Truncated Icosidodecahedron   &$30^{12}$   &$(0,10), (0,2)^{10}$&$40^{72}$  &$(0,2)^{10}, (0,10)^{2}$\\
Snub Cube                     &$30_{3,0}^4$&$(4,4)^{3}$         &$40^{24}$  &$(0,2)^{4}, (2,2)^{4}, (0,16)$\\
Snub Dodecahedron             &$50_{5,0}^6$&$(4,4)^{5}$         &$200^{12}$ &$(12,12)^{5}, (0,80)$\\
\hline
\hline
\end{tabular}
\end{center}
}
\caption{$z$-structure of prisms on Platonic and Archimedean solids}
\label{PrismPlatonicArchimedean}
\end{table}

\begin{table}
\begin{center}
\begin{tabular}{||c|c|c|c||}
\hline
\hline
dimension& half-$d$-cube                      & $z$-vector          &int. vectors\\
\hline
$2$&$\frac{1}{2}\gamma_3=\alpha_3$  &$4^3$    &$(1,1)^2$\\
$3$&$\frac{1}{2}\gamma_4=\beta_4$   &$8^{24}$ &$(0,2)^{4}$\\
$4$&$\frac{1}{2}\gamma_5=Med(\beta_{5})$ &$12^{240}$ &$(0,1)^{8}, (0,2)^{2}$\\
$5$&$\frac{1}{2}\gamma_{6}$         &$32^{1440}$& $(0,2)^{4}, (0,3)^{8}$\\
$6$&$\frac{1}{2}\gamma_{7}$         &$120^{6720}$&$(0,3)^{24}, (0,12)^{4}$\\
$7$&$\frac{1}{2}\gamma_{8}$         &$36^{430080}$&$(0,2)^{12}, (0,4)^{3}$\\
$8$&$\frac{1}{2}\gamma_{9}$         &$84^{3870720}$&$(0,4)^{6}, (0,5)^{12}$\\
$9$&$\frac{1}{2}\gamma_{10}$        &$192^{38707200}$& $(0,5)^{24}, (0,12)^{6}$\\
$10$&$\frac{1}{2}\gamma_{11}$        &$216^{851558400}$&$(0,3)^{48}, (0,18)^{4}$\\
$11$&$\frac{1}{2}\gamma_{12}$        &$160^{30656102400}$    &$(0,6)^{8}, (0,7)^{16}$\\
$12$&$\frac{1}{2}\gamma_{13}$        &$880^{159411732480}$    &$(0,7)^{80}, (0,40)^{8}$\\
\hline
\hline
\end{tabular}
\end{center}
\caption{$z$-structure of half-$d$-cubes for $d\leq 13$}
\label{HalfCubeResult}
\end{table}

\begin{table}
\begin{center}
\begin{tabular}{||c|c|c||}
\hline
\hline
Schl\"afli symbol of $P$                      &$|Aut(P)|$     &$z$-vector\\
\hline
$\{\frac{5}{2}, 5\}$                          &$120$     &$6^{10}$\\
$\{\frac{5}{2}, 3\}$                          &$120$     &$10^6$\\
\hline
$\{\frac{5}{2}, 5, 3\}$                       &$14400$   &$20^{360}$\\
$\{5, \frac{5}{2}, 5\}$                       &$14400$   &$15^{480}$\\
$\{\frac{5}{2}, 3, 5\}$                       &$14400$   &$12^{600}$\\
$\{\frac{5}{2}, 5, \frac{5}{2}\}$             &$14400$   &$15^{480}$\\
$\{3, \frac{5}{2}, 5\}$                       &$14400$   &$20^{360}$\\
$\{\frac{5}{2}, 3, 3\}$                       &$14400$   &$30^{240}$\\
\hline
\hline
\end{tabular}
\end{center}
\caption{$z$-structure of non-convex regular $3$- and $4$-polytopes (adapted from pages 292 and 294 of \cite{Cox73})}
\label{CoxeterTable}
\end{table}

The {\em pyramid} operation $Pyr({\cal K})$ (respectively,
{\em bipyramid} operation $BPyr({\cal K})$) on a $d$-dimensional complex
${\cal K}$ is the $(d+1)$-dimensional complex 
obtained by adding one (respectively, two) new vertices, connected to all vertices of the original complex.


All $92$ {\em Johnson solids}, i.e. regular-faced $3$-polytopes
were found in \cite{Jo1}.
All regular-faced, but not semiregular, $d$-polytopes, $d\geq 4$
are known also (\cite{BlBl2}). This list consist of two infinite families
of $d$-polytopes ($Pyr(\beta_{d-1})$ and $BPyr(\alpha_{d-1})$),
three particular $4$-polytopes ($Pyr(Ico)$, $BPyr(Ico)$ and 
the union of $0_{21}+Pyr(\beta_{3})$, where $\beta_{3}$ is a facet 
of $0_{21}$) and, finally, any $4$-polytope (except of snub $24$-cell),
arising from $600$-cell by the following special cut of vertices.
If $E$ is a subset of the $120$ vertices of $600$-cell, such that any
two vertices in $E$ are not adjacent, then this polytope is
the convex hull of all vertices of $600$-cell, except those in $E$.

Conway \cite{Conway} enumerated all {\em Archimedean $4$-polytopes},
i.e. those having a vertex-transitive group of symmetry and whose cells
are regular or Archimedean polyhedra and prisms or antiprisms with
regular faces. The list consists of:
\begin{enumerate}
\item $45$ polytopes obtained by Wythoff's kaleidoscope construction 
from regular $4$-polytopes (see Table \ref{WythoffResult} and, more
generally, Section \ref{WythoffSection});

\item $17$ prisms on Platonic, other than Cube, and Archimedean solids (see Table \ref{PrismPlatonicArchimedean});

\item prisms on $APrism_m$ for any $m>3$ (see Conjecture \ref{ConjecturePrismOfAntiprism});

\item a doubly infinite set of $4$-polytopes, which are direct
products $C_p\times C_q$ of two regular polygons (if one of polygons is a
square, then one gets prisms on $Prism_m$) (see Conjecture \ref{ProductOfPolygons});

\item the snub 24-cell $s(3,4,3)$ (see Table \ref{RegularResult});

\item a $4$-polytope, called in \cite{Conway} {\em Grand Antiprism};
it has $100$ vertices (all from $600$-cell), $300$ cells $\alpha_3$ and $20$
cells $APrism_5$ 
(those antiprisms form two interlocking tubes).
\end{enumerate}

\begin{remark}\label{GrandAntiprism}
The Grand Antiprism has $z$-vector $30^{20}, 50^{40}, 90^{20}$. The corresponding intersection vectors are $(0,1)^{10}, (0,2)^{10}$ and $(0,1)^{10}, (0,2)^{20}$ and $(0,1)^{10}$, $(0,2)^{10}$, $(0,4)^{5}$, $(4,4)^{5}$.
\end{remark}


\begin{remark}\label{WebPageInformation}
Complete information on z-structure of $92$-Johnson polyhedra is available from \cite{WebPageRegular}. We found $25$ $z$-uniform ones.
\end{remark}

\begin{remark}\label{WebSpecialCut}
The number of polytopes, obtained by special cuts, is unknown but it is finite.
By special cutting with $1$, \dots, $7$ vertices, one obtains, respectively,
$1$, $7$, $436$, $4776$, $45775$, $334380$ polytopes. We expect that for
$24$ vertices, there is only one possible special cut, which
yields semiregular snub $24$-cell. For more than $25$ vertices, there
is no special cut possible (i.e. the skeleton of $600$-cell has 
independence number $24$, see page 82 of \cite{martini}).
Due to the difficulty of the computation and very large
size of data, we computed the $z$-structure of
special cuts of $600$-cell only up to $3$ vertices. Results are
available from \cite{WebPageRegular}.
\end{remark}

\begin{remark}
In Table \ref{CoxeterTable} note that:

(i) Amongst those eight polytopes only $\{5, \frac{5}{2}, 5\}$ and $\{\frac{5}{2}, 5, \frac{5}{2}\}$ are self-dual.

(ii) In the case of Great Stellated Dodecahedron $\{\frac{5}{2}, 3\}$,
the item $h$ in Table 1 on page 292 of \cite{Cox73} (corresponding to the length of a zigzag) was $\frac{10}{3}$, while in Table \ref{CoxeterTable}, we put the value $10$. In fact, our notion is combinatorial, while Coxeter define Petrie polygon as a skew polygon (see Figure 6.1A on page 93 of \cite{Cox73}).
\end{remark}

\section{General results on $z$-structure of some generalizations of regular polytopes}

\begin{proposition}\label{PropositionInfiniteRegular}
For infinite series of regular polytopes we have:

(i) $z(\alpha_d)=(d+1)^{d!/2}$ with $Int=(0,1)^{d+1}$ for $d\geq 4$ and $(1,1)^2$ for $d=3$.

(ii) $z(\beta_d)=(2d)^{2^{d-2}(d-1)!}$ with $Int=(0,2)^d$.
\end{proposition}
\proof Both polytopes are regular polytopes. Therefore, they are $z$-uniform.
In order to know the length of a zigzag, one needs to compute the successive images of a flag under $T=\sigma_d\dots\sigma_2\sigma_1$.

Denote by $\{0,\dots, d\}$ the vertices of $\alpha_d$. It is easy to see that the image of the flag $f=(\{0\}$, $\{0,1\}$,\dots, $\{0,\dots, d-1\})$ is $(\{1\}$,$\{1,2\}$,\dots, $\{1,\dots, d\})$, i.e. it is the image of $f$ under a cycle of length $d+1$. Therefore, its length is $d+1$ and there is no self-intersection; hence, the $z$-vector is as in (i).
Also, one can check that two different zigzags intersect at most once if $d\geq 4$. Hence, the intersection vector is $(0,1)^{d+1}$. The case $d=3$ is trivial.

Denote by $\pm e_i$ with $1\leq i\leq d$ the vertices of $\beta_d$. It is easy to see that the image of the flag $f=(\{e_1\}$,$\{e_1,e_2\}$,\dots, $\{e_1,\dots e_d\})$ is the flag $f'=(\{e_2\}$,$\{e_2,e_3\}$,\dots, $\{e_2,\dots, e_d\}$,$\{-e_1,e_2,\dots,e_d\})$. Denote by $\phi$ the composition of the cycle $(1,\dots, d)$ on the coordinates with the symmetry $(x_1,\dots, x_d)\mapsto (-x_1,x_2,\dots, x_d)$. The order of $\phi$ is $2d$ and $\phi(f)=f'$. Therefore, all zigzags have length $2d$ and there is no self-intersection.
If two zigzags are intersecting, then they, moreover, intersect twice, since $\phi^d=-Id$ and one gets $Int=(0,2)^d$. \qed

In Table \ref{HalfCubeResult} are given $z$-structure of half-$d$-cubes
for $d \le 13$; note that the length of any zigzag there divides $2(d-2)$.

\begin{proposition}
For half-$d$-cube it holds:

(i) There are $d!2^{d-1}(d-2)$ flags, forming one orbit for $d=3,4$ and $d-2$ orbits for $d\geq 5$.

(ii) It is $z$-uniform.
\end{proposition}
\proof Let us write the set of vertices of $\frac{1}{2}\gamma_d$ as $\{S\subset \{1, \dots, d\}\mbox{~with~}|S|\mbox{~even}\}$. One has $\frac{1}{2}\gamma_3=\alpha_3$ and $\frac{1}{2}\gamma_4=\beta_4$, which are regular polytopes and whose structure is known. Therefore, one can assume $d\geq 5$.
The list of facets of $\frac{1}{2}\gamma_d$ consists of:
\begin{enumerate}
\item $2d$ facets $x_i=0$ and $x_i=1$
(those facets are incident to $2^{d-2}$ vertices of 
$\frac{1}{2}\gamma_d$, which form a polytope $\frac{1}{2}\gamma_{d-1}$).
\item $2^{d-1}$ simplex facets generated by vertices $\{S_1, \dots, S_{d}\}$ with $|S_i\Delta S_j|=2$ if $i\not=j$.
\end{enumerate}
From the above list of facets, one can easily deduce the list of $i$-faces of $\frac{1}{2}\gamma_d$; they are:
\begin{enumerate}
\item all $\frac{1}{2}\gamma_i$ with $4\leq i\leq d-1$ and 
\item all $k$-sets $\{S_1, \dots, S_k\}$ with $|S_i\Delta S_j|=2$ if $i\not= j$.
\end{enumerate}
The first kind of faces is obtained by intersecting hyperplanes $x_l=0,1$, while the second is obtained by taking any subset of a simplex face of $\frac{1}{2}\gamma_d$.
The symmetry group of $\frac{1}{2}\gamma_d$ has size $2^{d-1}d!$. It is generated by permutations of $d$ coordinates and operation $S\mapsto S_0\Delta S$ for a fixed $S_0\in \frac{1}{2}\gamma_d$.
There is one orbit of $k$-dimensional faces if $k\leq 2$ and two orbits, otherwise.

Take a flag $F_0\subset F_1\subset \dots\subset F_{d-1}$. If $F_i$ is a simplex face, then all faces, contained in it, are also simplexes. Therefore, the orbit, to which a flag belongs, is determined by the highest index $i$, for which it is still a simplex. Since $2\leq i\leq d-1$, this makes $d-2$ orbits. This yields (ii), since the stabilizer of a flag is trivial.

Let us denote by $O_i$ with $2\leq i\leq d-1$, the orbit formed by all flags, whose highest index is $i$. One has $\sigma_4(O_2)\subset O_3$ and $\sigma_k(O_2)\subset O_2$ for $k\not=2$. If $i=d-1$, then $\sigma_d(O_{d-1})\subset O_{d-2}$, while $\sigma_k(O_{d-1})\subset O_{d-1}$ if $k\not= d$. If $2<i<d$, then one has $\sigma_{i+2}(O_i)\subset O_{i+1}$ and $\sigma_{i+1}(O_i)\subset O_{i-1}$; for other $k$, one has $\sigma_k(O_i)\subset O_i$.

Recalling $T=\sigma_d\sigma_{d-1}\dots\sigma_1$, one obtains $T(O_i)\subset O_{i-1}$ if $i>2$ and $T(O_2)\subset O_{d-1}$. Therefore, all orbits of flags are touched by any zigzag of $\frac{1}{2}\gamma_d$. This proves $z$-uniformity. \qed


\begin{proposition}
For $Pyr(\beta_{d-1})$, it holds:

\begin{enumerate}
\item[(i)] there are $(d+1)(d-1)!2^{d-2}$ flags partitioned into $d+1$ orbits.
\item[(ii)] it is $z$-uniform.
\end{enumerate}

\end{proposition}
\proof Denote by $v$ the vertex, on which we do the pyramid construction.
Take a flag $(F_1, \dots, F_{d})$ of $Pyr(\beta_{d-1})$. The sequence of faces $(F_1\cap \beta_{d-1}, \dots, F_{d}\cap \beta_{d-1})$ 
can not be a flag for three possible reasons:
\begin{enumerate}
\item $F_1\cap \beta_{d-1}=\emptyset$, it means that $F_1=\{v\}$.
\item $F_i\cap \beta_{d-1}=F_{i+1}\cap \beta_{d-1}$, it means that $F_{i+1}=conv(F_i, v)$.
\item $F_{d}\cap\beta_{d-1}=\beta_{d-1}$, it means that $F_{d}=\beta_{d-1}$.
\end{enumerate}
This implies, since $\beta_{d-1}$ is regular, that $Pyr(\beta_{d-1})$ has the following orbits of flags:
\begin{enumerate}
\item $O_i$, with $1\leq i\leq d$, being the orbit of flags of $Pyr(\beta_{d-1})$, whose first face containing $v$ is in position $i$;
\item the orbit $O_{d+1}$ of flags obtained by adding $\beta_{d-1}$ to a flag of $\beta_{d-1}$.
\end{enumerate}

The operator $\sigma_i$ with $1\leq i\leq d$, which acts on the flag $(F_1, \dots, F_{d})$ by exchanging the term $F_{i}$, acts on the orbit by permuting the orbits $O_i$ and $O_{i+1}$ and leaving the others preserved. Hence, the product $T$ acts on the set of orbits $O_i$ as the cycle $(1,2,\dots, d+1)$. So, $Pyr(\beta_{d-1})$ is $z$-uniform. \qed

\begin{conjecture}\label{ConjectureBetaAlpha}
(i) For $z$-structure of $Pyr(\beta_{d-1})$ it holds:

(i.1) $z$-vector is:
\begin{equation*}
\left\lbrace\begin{array}{rcl}
(d^2-1)^{(d-2)!2^{d-2}}&\mbox{~for~}&d\mbox{~even},\\
{2(d^2-1)}_{2d-2,0}^{(d-2)!2^{d-3}}&\mbox{~for~}&d\mbox{~odd~and~}d>3,\\
16_{8,8}&\mbox{~for~}&d=3.
\end{array}\right.
\end{equation*}

(i.2) Intersection vectors are:
\begin{equation*}
\left\lbrace\begin{array}{rcl}
(0,d-1)^{d-1}, (0,2d-2)&\mbox{for~}&d\mbox{~even~and~}d\geq 4,\\
(0,2d-2)^{d-1}         &\mbox{for~}&d\mbox{~odd}.
\end{array}\right.
\end{equation*}

(ii) For $z$-structure of $BPyr(\alpha_{d-1})$ it holds:

(ii.1) $z$-vector is:
\begin{equation*}
\left\lbrace\begin{array}{rcl}
(d^2)^{(d-1)!}                    &\mbox{~for~}&d\mbox{~even~and~}d\geq 4,\\
({2d^2}_{2d,0})^{\frac{(d-1)!}{2}}  &\mbox{~for~}&d\mbox{~odd~and~}d>3,\\
(18_{6,3})                        &\mbox{~for~}&d=3.
\end{array}\right.
\end{equation*}

(ii.2) Intersection vectors are:
\begin{equation*}
\left\lbrace\begin{array}{rcl}
(0,2d), (0,d-2)^{d}    &\mbox{~for~}&d\mbox{~even~and~}d>4,\\
(0,2d-4)^d              &\mbox{~for~}&d\mbox{~odd~and~}d>3,\\
(0,8),(2,2)^2&\mbox{~for~}&d=4,\\
\end{array}\right.
\end{equation*}

\end{conjecture}
Clearly, $Pyr(\beta_2)$ and $BPyr(\alpha_2)$ are just square pyramid and 
dual $Prism_3$, respectively. Above conjecture was checked for $n\leq 10$.

\begin{conjecture}\label{ProductOfPolygons}
Let $t$ denote $gcd(p,q)$ and $s$ denote $\frac{pq}{t^2}$. Then for
$z$-structure of the direct product $C_p \times C_q$ holds:

(i) If $p,q$ are both even, then

$z=(2ts)^{6t}$ with $Int=(0,2s)^t$ for all zigzags.

(ii) If exactly one of $p,q$ is odd, then

$z=(2ts)^{6t}$ with $Int=(0,s)^{2t}$ for $4t$ zigzags and $Int=(s,s)^t$ for
the remaining $2t$ zigzags.

(iii) If $p,q$ are both odd, then

$z=(2ts)^{2t},(4ts)^{2t}$ with $Int=(s,s)^t$ for zigzags of length $2ts$
and $Int=(2s,2s)^t$ for zigzags of length $4ts$.
\end{conjecture}
The above conjecture was checked for $p,q\leq 15$.

For any zigzag of $Prism(P)$ with $z(P)=a^b$ and $P$ 
being Platonic or Archimedean $3$-polytope, one has
$z=(\frac{4a}{gcd(a,3)})^{2gcd(a,3)b}$.
In general, $z=(\frac{da}{gcd(a,d-1)})^{2gcd(a,d-1)b}$.
Cube is not included in Table \ref{PrismPlatonicArchimedean},
because the prism on it is just $\gamma_4$.
Above relation works also for prisms on antiprisms.


\begin{conjecture}\label{ConjecturePrismOfAntiprism}
For $z$-structure of prism on $APrism_m$ it holds:

$z=(\frac{8m}{gcd(m,3)})^{8gcd(m,3)}$ with $Int=(\frac{2m}{3})^4$ if $gcd(m,3)=3$ and, otherwise,
two zigzags have $Int=(0,2m)^4$,
two zigzags have $Int=(0,2m),(2m,4m)$ and
four zigzags have $Int=(0,2m)^2,(0,4m)$.
\end{conjecture}
The above conjecture was checked for $m\leq 15$.

Denote by $I(Z_1, Z_2)=(n_I, n_{II})$ the pair of intersection numbers 
between two zigzags, $Z_1$ and $Z_2$, corresponding to
intersections of type I and II.
Given a map $f$ acting on a complex ${\cal K}$ without any fixed face, the
{\em folded complex} $\tilde{{\cal K}}$ is defined as the quotient space of
${\cal K}$ under $f$; it is not always a lattice.

\begin{proposition}\label{FixedPointSymmetry}
Let ${\cal K}$ be a complex and $f$ a fixed-point free involution on ${\cal K}$; then one has:

(i) For any zigzag $Z$ of ${\cal K}$, such that $f(Z)=Z$, the length and the 
signature of its image $\tilde{Z}$ in $\tilde{K}$ are the half of the length
and the signature, respectively, of $Z$.

(ii) If $Z_2=f(Z_1)$ with $Z_2\not= Z_1$, then we put compatible orientation on $Z_1$ and $Z_2$. The zigzags $Z_1$ and $Z_2$ are mapped to a zigzag $\tilde{Z}$ of $\tilde{{\cal K}}$ with its signature being equal to the signature of $Z_1$ plus $\frac{1}{2}I(Z_1, Z_2)$.

Concerning intersection vectors, one has:

(i) Two zigzags of ${\cal K}$, which are invariant under $f$, are mapped to zigzags of $\tilde{{\cal K}}$ with halved intersection.

(ii) Take an invariant zigzag $Z$ of ${\cal K}$ and $Z_2=f(Z_1)$ two equivalent zigzags of ${\cal K}$. They are mapped to $\tilde{Z}$ and $\tilde{Z'}$ and one has  $I(\tilde{Z}, \tilde{Z'})=I(Z, Z_1)$.


(iii) Take two pairs $(Z_1, Z'_1)$ and $(Z_2, Z'_2)$ with $Z'_i=f(Z_i)$. They are mapped to $\tilde{Z_1}$ and $\tilde{Z_2}$ and their intersection $I(\tilde{Z_1}, \tilde{Z_2})$ is equal to $I(Z_1, Z_2)+I(Z_1, Z'_2)$.


\end{proposition}

For example, Petersen graph, embedded on the projective plane, is a folding of the Dodecahedron by central inversion. Another example is a map on torus, which is folded onto the Klein bottle.

%


The {\em folded cube} $\Box_d$ is obtained from $d$-cube by
{\em folding}, i.e. by identifying opposite faces of 
$\gamma_d$. Obtained complex is $(d-1)$-dimensional, like $\gamma_d$, 
but it is not a lattice, which imply that this complex
does not admit a realization as {\em polyhedral complex}.

\begin{proposition}
For $\Box_d$ one has $z=d^{2^{d-2}(d-1)!}$ with $Int=(0,1)^d$.
\end{proposition}
\proof Every zigzag of $\beta_d$ corresponds to a zigzag of $(\beta_d)^*=\gamma_d$; hence, by the proof of Proposition \ref{PropositionInfiniteRegular}, the zigzags of $\gamma_d$ are centrally symmetric. By applying Theorem \ref{FixedPointSymmetry}, one obtains $z=d^{2^{d-2}(d-1)!}$. Furthermore, one can prove easily that zigzags of $\gamma_d$ have $Int=(0,2)^d$; hence, the intersection vector of $\Box_d$ is $(0,1)^d$. \qed

A $(d-1)$-dimensional complex ${\cal K}$ is said to be {\em of type $\{3,4\}$} if every $(d-2)$-dimensional face is contained in $3$ or $4$ faces of dimension $d-1$. Those simplicial complexes are classified in terms of partitions: given such a simplicial complex, there exist a partition $(P_1, \dots, P_t)$ of $\{1, \dots, d\}$, such that ${\cal K}^*$ is isomorphic to $\Delta_1\times \Delta_2\times\dots\times\Delta_t$ with $\Delta_i$ being the simplex of dimension $|P_i|$; see \cite{DDS} for details.

\begin{conjecture}

(i) A simplicial complex of type $\{3,4\}$ is not $z$-uniform if and only if the sizes of parts in the corresponding partition are either $(\frac{d}{2},\frac{d}{2})$, or all even (except simplex).

In non-$z$-uniform case, $gcd(l_1,l_2)=min(l_1,l_2)$ for any two lengths of zigzags.

(ii) In special case $\{1,\dots,\frac{d}{2}\},\{\frac{d}{2}+1,\dots,d\}$ one has $max(l_i)=\frac{d(d+2)}{2}$ and $min(l_i)=d+2$. In other extreme case $\{1,2\},\dots,\{d-1,d\}$ one has $max(l_i)=3d$ and $min(l_i)=\frac{3d}{2}$.

(iii) For partition $\{1\},\{2,\dots, d\}$ the simplicial complex of type $\{3,4\}$ is, in fact, $BPyr(\alpha_{d-1})$.
%
%

(iv) For partition $\{1\},\{2\},\dots,\{d-2\},\{d-1,d\}$ the simplicial complex of type $\{3,4\}$ has the following $z$-structure:

(iv.1) $\lfloor \frac{d}{2}\rfloor$ orbits, each zigzag has length $6d$ and intersection vector $(d,0)^6$ (intersection vectors are $(12,6)$ for $d=3$ and $(0,8),(2,2)^2$ for $d=4$).

(iv.2) For odd $d$, all orbits have $2^{d-3}(d-2)!$ zigzags. For even $d$, one orbit has size $2^{d-4}(d-2)!$ and $\frac{d-2}{2}$ orbits have size $2^{d-3}(d-2)!$.
\end{conjecture}
The above conjecture was checked up to $d=8$.



\begin{table}
\begin{center}
{\scriptsize
\begin{tabular}{||c|c|c||}
\hline
\hline
Wythoff $4$-polytope                            & $z$-vector          &intersection vectors\\
\hline
$\alpha_4$=$\alpha_4(\{0\})$=$\alpha_4(\{3\})$          &$5^{12}$             &$(0,1)^5$\\
$\alpha_4(\{0,1\})$=$\alpha_4(\{2,3\})$       &$20^{12}$            &$(0,4)^5$\\
$\alpha_4(\{0,1,2\})$=$\alpha_4(\{1,2,3\})$   &$20^{36}$            &$(0,1)^5, (0,3)^{5}$ or $(1,3)^{5}$\\
$\alpha_4(\{0,1,2,3\})$                  &$20^{72}$            &$(0,2)^{10}$\\
$\alpha_4(\{0,1,3\})$=$\alpha_4(\{0,2,3\})$   &$48^{20}$            &$(0,5)^{6}, (3,15)$\\
$\alpha_4(\{0,2\})$=$\alpha_4(\{1,3\})$       &$45^{12}$            &$(4,5)^{5}$\\
$\alpha_4(\{0,3\})$                      &$10^{12}, 30^{12}$   &$(0,2)^5$ or $(0,6)^5$\\
$0_{21}$=$\alpha_4(\{1\})$=$\alpha_4(\{2\})$  &$15^{12}$            &$(1,2)^5$\\
$\alpha_4(\{1,2\})$                      &$10^{12}, 20^{12}$   &$(0,2)^5$ or $(0,4)^5$\\
$\beta_4$=$\beta_4(\{0\})$              &$8^{24}$             &$(0,2)^4$\\
$\beta_4(\{0,1\})$                      &$16^{48}$            &$(0,1)^8, (0,4)^2$\\
$\beta_4(\{0,1,2\})$=$24$-cell$(\{0,1\})$=$24$-cell$(\{2,3\})$&$24^{96}$&$(0,2)^9, (0,6)$\\
$\beta_4(\{0,1,2,3\})$                  &$32^{144}$           &$(0,2)^{16}$\\
$\beta_4(\{0,1,3\})$                    &$64^{48}$            &$(0,2)^{2}, (0,4)^{4}, (0,6)^{4}, (4,6)^{2}$\\
$\beta_4(\{0,2\})$=$24$-cell$(\{1\})$=$24$-cell$(\{2\})$&$18^{96}$      &$(0,2)^{9}$\\
$\beta_4(\{0,2,3\})$                    &$64^{48}$            &$(0,2)^{2}, (0,6)^{10}$\\
$\beta_4(\{0,3\})$                      &$16^{24}, 48^{24}$   &$(0,2)^{8}$ or $(0,6)^{8}$\\
$24$-cell=$\beta_4(\{1\})$=$24$-cell$(\{0\})$=$24$-cell$(\{3\})$&$12^{48}$        &$(0,2)^{6}$\\
$\beta_4(\{1,2\})$                      &$16^{24}, 24^{32}$   &$(0,2)^{8}$ or $(0,2)^{6}, (0,4)^{3}$\\
$\beta_4(\{1,2,3\})$                    &$32^{72}$            &$(0,2)^{4}, (2,4)^{4}$ or $(0,2)^{8}, (0,4)^{4}$\\
$\beta_4(\{1,3\})$                      &$36^{48}$            &$(0,2)^{8}, (0,4)^{3}, (0,8)$\\
$\beta_4(\{2\})$                        &$24^{24}$            &$(0,2)^{4}, (0,4)^{4}$\\
$\beta_4(\{2,3\})$                      &$32^{24}$            &$(0,2)^{4}, (0,6)^{4}$\\
$\gamma_4$=$\beta_4(\{3\})$             &$8^{24}$             &$(0,2)^{4}$\\
$24$-cell$(\{0,1,2\})$=$24$-cell$(\{1,2,3\})$   &$48^{144}$           &$(0,2)^{12}, (0,4)^{6}$ or $(0,2)^{6}, (2,4)^{6}$\\
$24$-cell$(\{0,1,2,3\})$                  &$48^{288}$           &$(0,2)^{24}$\\
$24$-cell$(\{0,1,3\})$=$24$-cell$(\{0,2,3\})$     &$96^{96}$            &$(0,4)^{6}, (0,6)^{7}, (4,6)^{3}$\\
$24$-cell$(\{0,2\})$=$24$-cell$(\{1,3\})$       &$54^{96}$            &$(0,12), (0,2)^{3}, (0,4)^{6}, (0,6)^{2}$\\
$24$-cell$(\{0,3\})$                      &$24^{192}$           &$(0,2)^{12}$\\
$24$-cell$(\{1,2\})$                      &$24^{48}, 48^{48}$   &$(0,2)^{12}$ or $(0,4)^{12}$\\
$600$-cell=$600$-cell$(\{0\})$             &$30^{240}$           &$(0,2)^{15}$\\
$600$-cell$(\{0,1\})$                      &$48^{600}$           &$(0,2)^{24}$\\
$600$-cell$(\{0,1,2\})$                    &$80^{1080}$          &$(0,2)^{40}$\\
$600$-cell$(\{0,1,2,3\})$                  &$120^{1440}$         &$(0,2)^{60}$\\
$600$-cell$(\{0,1,3\})$                    &$320^{360}$          &$(0,4)^{20}, (2,4)^{10}, (6,12)^{10}$\\
$600$-cell$(\{0,2\})$                      &$135^{480}$          &$(0,2)^{15}, (0,3)^{15}, (0,4)^{15}$\\
$600$-cell$(\{0,2,3\})$                    &$192^{600}$          &$(0,2)^{30}, (0,4)^{12}, (2,12)^{6}$\\
$600$-cell$(\{0,3\})$                      &$60^{960}$           &$(0,2)^{30}$\\
octicosahedric polytope=$600$-cell$(\{1\})$           &$45^{480}$           &$(0,1)^{15}, (0,2)^{15}$\\
$600$-cell$(\{1,2\})$                      &$60^{240}, 80^{360}$ &$(0,2)^{30}$ or $(0,2)^{20}, (0,4)^{10}$\\
$600$-cell$(\{1,2,3\})$                    &$120^{720}$          &$(0,2)^{15}, (2,4)^{15}$\\
                                    &                     &or $(0,2)^{30}, (0,4)^{15}$\\
$600$-cell$(\{1,3\})$                      &$108^{600}$          &$(0,2)^{12}, (0,4)^{6}, (2,8)^{6}$\\
$600$-cell$(\{2\})$                        &$90^{240}$           &$(0,2)^{15}, (0,4)^{15}$\\
$600$-cell$(\{2,3\})$                      &$120^{240}$          &$(0,2)^{15}, (0,6)^{15}$\\
$120$-cell=$600$-cell$(\{3\})$             &$30^{240}$           &$(0,2)^{15}$\\
\hline
\hline
\end{tabular}
}
\end{center}
\caption{$z$-structure of Wythoff Archimedean $4$-polytopes}
\label{WythoffResult}
\end{table}

\section{Wythoff kaleidoscope construction}\label{WythoffSection}
Wythoff construction is defined for any $d$-dimensional complex ${\cal K}$ and non-empty subset $V$ of $\{0,\dots, d\}$. It was introduced in \cite{wythoff} and \cite{Cox35}.

The set of all {\em partial flags} $(f_{i_0}, \dots, f_{i_m})$, with $f_{i_j}\subset f_{i_{j+1}}$ and $i_j\in V$, is the vertex-set of a complex, which we denote by ${\cal K}(V)$ and call {\em Wythoff construction with respect to the complex ${\cal K}$ and the set $V$}.

In general, one has 
${\cal K}(V)={\cal K}^*(d-V)$ with $d-V$ denoting the set of all $d-i$, 
$i\in V$.
If a complex ${\cal K}$ is self-dual, then one has ${\cal K}(V)={\cal K}(d-V)$.
One has, in general, ${\cal K}(\{0\})={\cal K}$, ${\cal K}(\{d\})={\cal K}^*$
and ${\cal K}(\{1\})=Med({\cal K})$.
Dual ${\cal K}(\{0,\dots,n\})$ is a simplicial $(d-1)$-complex called {\em order-complex} (\cite{RS}).

Easy to see that a general $d$-dimensional complex admits at most $2^{d+1}-1$ non-isomorph Wythoff constructions, while a self-dual $d$-dimensional complex admits at most $2^{d}+2^{\lceil \frac{d-1}{2}\rceil}-1$ such non-isomorph constructions.
Curiously, in the regular complexes considered, we obtain exactly $2^{d+1}-1$ and $2^{d}+2^{\lceil \frac{d-1}{2}\rceil}-1$ non-isomorph complexes.


If ${\cal K}$ is a $2$-dimensional complex,
then it is easy to see that ${\cal K}(V)$ with $V$=$\{0\}$, $\{0,1\}$, $\{0,1,2\}$, $\{0,2\}$, $\{1,2\}$, $\{1\}$ and $\{2\}$ correspond, respectively, to following maps: original map ${\cal M}$, truncated ${\cal M}$, truncated  $Med({\cal M})$, $Med(Med({\cal M}))$, truncated  ${\cal M}^*$, $Med({\cal M})$ and ${\cal M}^*$ (dual ${\cal M}$).

Call {\em Wythoff Archimedean} any Wythoff construction with 
respect to some regular $d$-polytope.
By applying the Wythoff construction to the three $3$-valent Platonic
solids (Tetrahedron, Cube and Dodecahedron) one obtains all Archimedean $3$-polytopes, except Snub Cube and Snub Dodecahedron; their $z$-structure is indicated in columns 2 and 3 of Table \ref{PrismPlatonicArchimedean}.
See in Table \ref{WythoffResult} the $z$-structure of Wythoff Archimedean $4$-polytopes.



\section{Lins triality}\label{LinsSection}

\begin{table}
\begin{center}
\begin{tabular}{||c|c||c|c|c|c||}
\hline
\hline
$(v,f,z)$   &$(a,b,c)$   &our notation          &notation    &notation      &Euler\\
$\downarrow$&$\downarrow$&                      &in \cite{L} &in \cite{AS02}&characteristic\\
\hline
$(v,f,z)$   &$(a,b,c)$   &${\cal M}$            &gem                   &${\cal M}$&$\chi({\cal M})$\\  
$(f,v,z)$   &$(c,b,a)$   &${\cal M}^*$          &dual gem              &${\cal M}^*$&$\chi({\cal M})$\\  
$(z,f,v)$   &$(a, b, ac)$&$phial({\cal M})$     &phial gem             &$p((p({\cal M}))^*)$&$\chi_p({\cal M})$\\  
$(f,z,v)$   &$(ac,b,a)$  &$(phial({\cal M}))^*$ &skew-dual gem         &$(p({\cal M}))^*$&$\chi_p({\cal M})$\\  
$(v,z,f)$   &$(ac,b,c)$  &$skew({\cal M})$      &skew gem              &$p(M)$&$\chi_s({\cal M})$\\  
$(z,v,f)$   &$(c,b,ac)$  &$(skew({\cal M}))^*$  &skew-phial gem        &$p({\cal M}^*)$&$\chi_s({\cal M})$\\  
\hline
\hline
\end{tabular}
\end{center}
\caption{Lins triality}
\label{LinsDual}
\end{table}

In the case of maps on surfaces, flags are triples $(v,e,f)$ with $v\in e\subset f$, where $v$, $e$ and $f$ are incident vertex, edge, and face, respectively.
Denote by $a$, $b$ and $c$ the three mappings $\sigma_1$, $\sigma_2$ and $\sigma_3$.
Vertex, edge and face are identified with the set of flags containing them; therefore, with orbits on flags of the groups $\langle b, c\rangle$, $\langle a, c\rangle$ and $\langle a, b\rangle$.
Zigzags were defined in Section \ref{GeneralCase} above as circuits of
flags $(f_i)_{1\leq i\leq l}$ with $f_{i+1}=cba f_i$. It is easy to see that this correspond to orbits of the group $\langle ac, b\rangle$.

Let $v=(v_i)_{i\geq 1}$, $p=(p_j)_{j\geq 1}$ and $z=(z_k)_{k\geq 1}$ are, respectively, $v$-, $p$- and $z$-vectors of a map.
Then the number $\sum_{i\geq 1}iv_i=\sum_{j\geq 1}jp_j=\sum_{k\geq 1}kz_k$ is the double of the number of edges.

One can reconstruct the map from the flag-set and the triple $(a,b,c)$ of operations acting on it by using the representation of vertices, edges and faces as orbits.
The only restriction, that applies to $a$, $b$ and $c$ is $a^2=b^2=c^2=(ac)^2=1$.
If one changes the triple $(a,b,c)$ to $(c,b,a)$, then the map is changed to its dual.

%


%
%

Other operations were introduced in \cite{L}: mapping $(a,b,c)$ to $(a,b,ac)$ or $(ac,b,c)$ produce the maps called $phial({\cal M})$ and $skew({\cal M})$. In \cite{JT87} it is proved that 
there is no other ``good'' notions of dualities for maps on surfaces than
the six ones given in Table \ref{LinsDual}.
The {\em skeleton} graph of a map (i.e. the graph of its vertices and edges) is connected. It is well-known that the dual graph of any connected map on a surface is connected also.
By using operation $phial$, we see that, moreover, the graph of zigzags $Z({\cal M})$ is connected.

The six operations depicted in Table \ref{LinsDual} form a group isomorphic to $Sym(3)$. In particular, each of operations $dual$, $skew$ and $phial$ is a reflexion.

\begin{figure}
\begin{center}
\begin{minipage}{6cm}
\centering
\epsfig{file=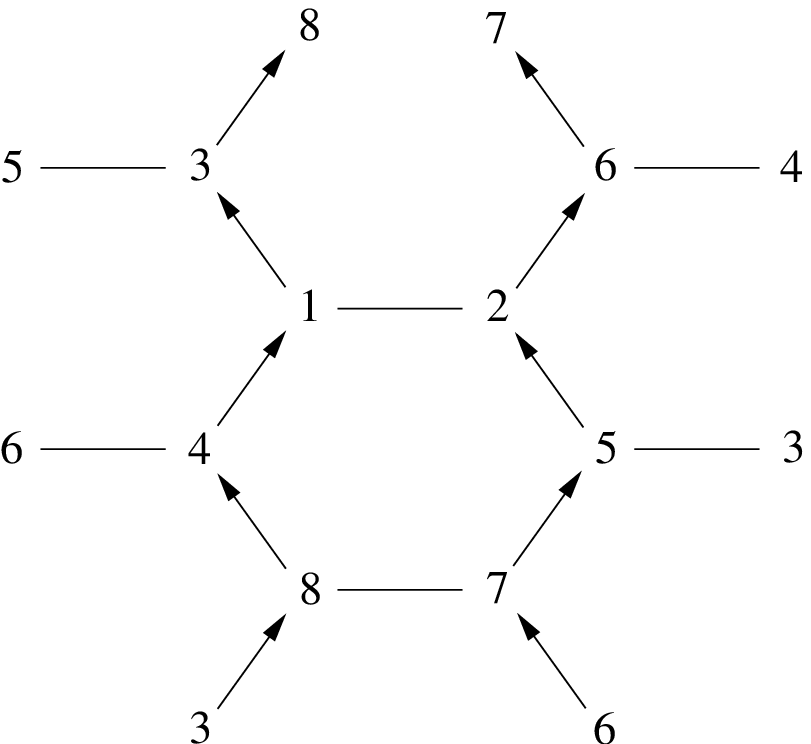, height=4cm}
\end{minipage}
\begin{minipage}{6cm}
\centering
\epsfig{file=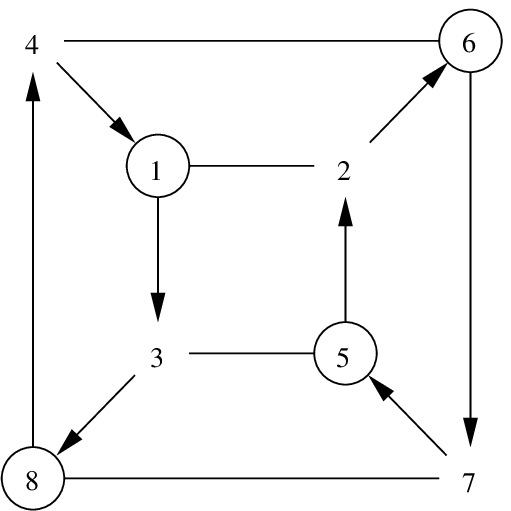, height=4cm}
\end{minipage}
\end{center}
\caption{Two representations of $skew(Cube)$: on torus and as a Cube with twisted (encircled) vertices.}
\label{ExampleOfTheCube}
\end{figure}


\begin{figure}
\begin{center}
\begin{minipage}{5cm}
\centering
\epsfig{file=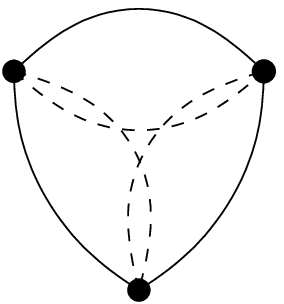, height=3cm}\par
$phial(Tetrahedron)$
\end{minipage}
\begin{minipage}{5cm}
\centering
\epsfig{file=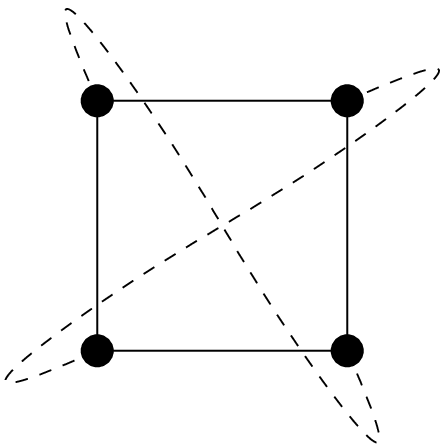, height=3cm}\par
$skew(Tetrahedron)$
\end{minipage}
\end{center}
\caption{Two maps on projective plane.}
\label{TwoMapsProjective}
\end{figure}

Denote the Euler characteristic of ${\cal M}$, $phial({\cal M})$ and $skew({\cal M})$ by $\chi({\cal M})$, $\chi_p({\cal M})$ and $\chi_s({\cal M})$, respectively.

\begin{conjecture}
(i) For Lins triality for $Prism_m$ it holds:

(i.1) $\chi_s(Prism_m)=gcd(m,4)-m$ and $skew(Prism_m)$ is oriented if and only if $m$ is even,

(i.2) $\chi_p(Prism_m)=2+gcd(m,4)-2m=\chi(Prism_m)+\chi_s(Prism_m)-m$ and $phial({\cal M})$ is non-oriented.

(ii) For Lins triality for $APrism_m$ it holds:

(ii.1) $\chi_s(APrism_m)=1+gcd(m,3)-2m$ and $skew({\cal M})$ is non-oriented,

(ii.2) $\chi_p(APrism_m)=3+gcd(m,3)-2m=\chi(APrism_m)+\chi_s(APrism_m)$~~~and $skew(APrism_m)$ is oriented.

\end{conjecture}
The above conjecture was checked up to $n=100$.

The $phial(Tetrahedron)$ is the complex obtained by taking the octahedron and identifying opposite points,
while $skew(Tetrahedron)$ is the complex obtained by taking Cube and identifying opposite points; see Figure \ref{TwoMapsProjective}.

The complex $skew(Cube)$ is a $3$-valent map on the torus with $8$ vertices and $4$ hexagonal faces (twisted construction); see Figure \ref{ExampleOfTheCube}.

A vertex of a graph, embedded in an orientable surface, is called {\em twisted} if the clockwise order of its adjacent vertices is the reversal, with respect of original clockwise order, given by the original embedding.


\begin{conjecture}
Let ${\cal M}$ be a map on an oriented surface, such that its
skeleton $G({\cal M})$ is bipartite, then $skew({\cal M})$
is a map on an oriented surface and $G(skew({\cal M}))=G({\cal M})$.
The orientation of surface induces, for each vertex $x$
of $G({\cal M})$, a cyclic order on vertices, to which $x$ is adjacent; 
then the maps ${\cal M}$ and $skew({\cal M})$ differ only by the twisting
of the vertices of one part of the bipartition of $G({\cal M})$.
\end{conjecture}

\begin{figure}
\begin{center}
\begin{minipage}[t]{5.5cm}
\centering
\epsfig{file=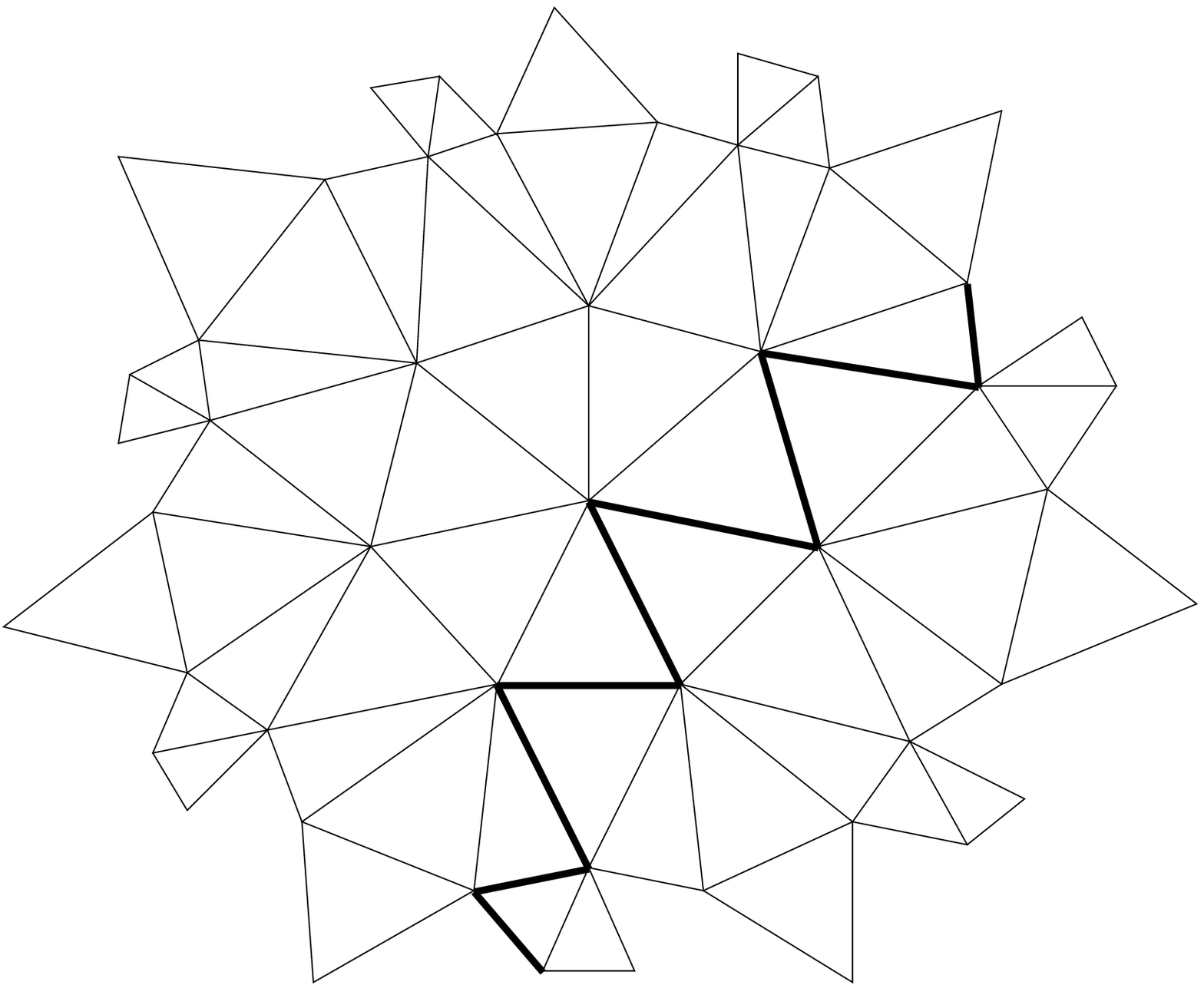, height=3.5cm}
\end{minipage}
\hspace{2cm}
\begin{minipage}[t]{5.5cm}
\centering
\epsfig{file=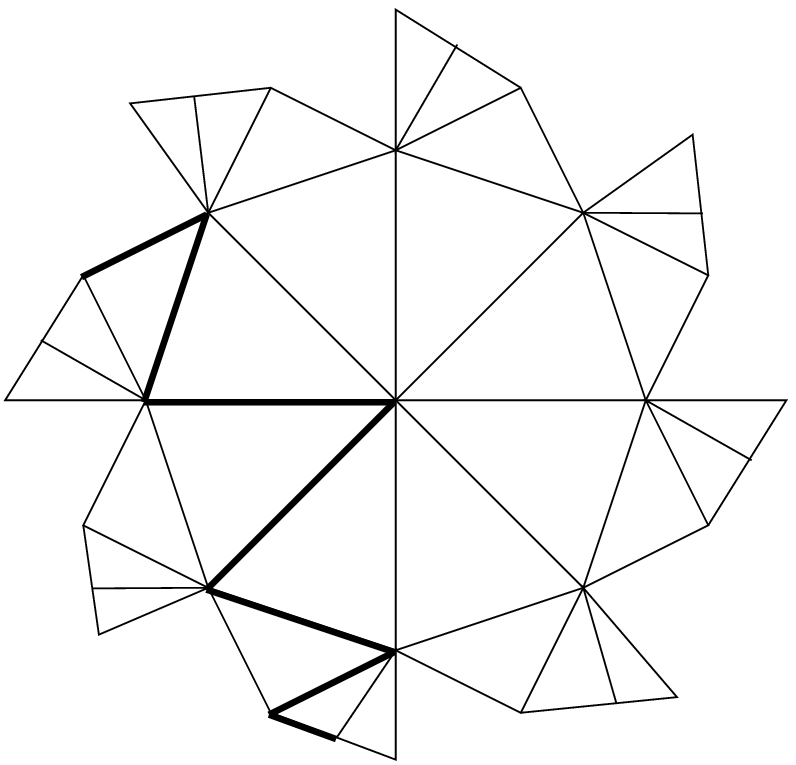, height=3.5cm}
\end{minipage}
\end{center}
\caption{Typical zigzags in dual Klein and dual Dyck maps}
\label{Dyck-Klein-Map}
\end{figure}

In particular, $Skew(Prism_{2k})$ is $Prism_{2k}$ with $k$ independent vertices being twisted, if above conjecture is true; we checked it for $k\leq 4$. Also, we checked above conjecture for two following cases:

(i) $(phial({\cal M}))^*$ (i.e. $skew\circ dual$) of $APrism_4$ is dual $APrism_4$ with exactly five independent vertices (i.e. a part of this bipartition of two parts of size $5$) being twisted.

(ii) $(phial(Cuboctahedron))^*$ is dual Cuboctahedron with exactly one part (eight $3$-valent vertices) of this bipartite graph (eight $3$-valent and six $4$-valent vertices) being twisted.

\begin{table}
\begin{center}
{\scriptsize
\begin{tabular}{||c|c|c|c|c|c||}
\hline
\hline
regular $3$-valent map   &Genus &Nr. vertices  &rotation group          &$z$-vector    &\multicolumn{1}{|c||}{$\QuotS{z(GC_{k,l})}{(k^2+kl+l^2)}$}\\
\hline
Dodecahedron $=\{5,3\}$   &$0$&$20$ &$A_5\simeq PSL(2,5)={}^{5}T$  &$10^6$         &$10^6$ or $6^{10}$ or $4^{15}$\\
dual Klein map $\{7,3\}$ &$3$&$56$ &$PSL(2,7)={}^{7}O$      &$8^{21}$       &$6^{28}$ or $8^{21}$\\
dual Dyck map $\{8,3\}$  &$3$&$32$ &${}^{4}O$ &$6^{16}$     &$6^{16}$ or $8^{12}$\\
$\{11,3\}$               &$26$&$220$&$PSL(2,11)={}^{11}I$    &$10^{66}$       &$6^{110}$ or $10^{66}$ or $12^{55}$\\
\hline
\hline
\end{tabular}
}
\end{center}
\caption{$z$-structure of some regular $3$-valent maps and of their Goldberg-Coxeter $GC_{k,l}$ construction (see \cite{goldbergcoxeter})}
\label{TableRegMaps}
\end{table}

In Table \ref{TableRegMaps} are presented several regular maps. 
The group of dual Dyck map is denoted by ${}^{4}O$,
because $O$ is a subgroup of index $4$ of it; by the same reason, this
group (of order 96) is called {\em tetrakisoctahedral}; it is generated by
two elements $R$, $S$ subject to the relations
$R^3=S^8=(RS)^2=(S^2R^{-1})^3=1$.
Now, $PSL(2,p)$ for $p=5,7,11$ are denoted by ${}^{5}T$, ${}^{7}O$,
${}^{11}I$, respectively, and called, respectively, {\em pentakistetrahedral},
{\em heptakioctohedral}, {\em undecakisicosahedral}.
A well-known result of Evariste Galois is that
they are the only ones amongst all $PSL(2,p)$, 
which act transitively on less than $p+1$ elements.


\end{document}